\documentclass[12pt,twoside]{article}
\usepackage{amsmath}

\date{}
\begin{document}

\centerline{}

\centerline{\Large{\bf Role of operator semirings in characterizing
}}

\centerline{\Large{\bf $\Gamma-$semirings in terms of fuzzy
subsets}}

\centerline{}

\centerline{\bf{Sujit Kumar Sardar}}

\centerline{Department of Mathematics,}

\centerline{Jadavpur University, Kolkata}

\centerline{E-mail: sksardarjumath@gmail.com}

\centerline{}

\centerline{\bf{Y. B. Jun}}

\centerline{Department of Mathematics Education}

\centerline{Gyeongsang National University}

\centerline{Chinju 660-701, Korea}

\centerline{skywine@gmail.com}

\centerline{}

\centerline{\bf{Sarbani Goswami}}

\centerline{Lady Brabourne College, Kolkata, W.B., India}

\centerline{E-mail: sarbani7$_{-}$goswami@yahoo.co.in}

\newtheorem{Theorem}{\quad Theorem}[section]

\newtheorem{Definition}[Theorem]{\quad Definition}

\newtheorem{Corollary}[Theorem]{\quad Corollary}

\newtheorem{Lemma}[Theorem]{\quad Lemma}

\newtheorem{Example}[Theorem]{\emph{Example}}

\newtheorem{Proposition}[Theorem]{Proposition}

\newtheorem{Remark}[Theorem]{\quad Remark}

\begin{abstract}
The operator semirings of a $\Gamma-$semiring have been brought into
use to study $\Gamma-$semiring in terms of fuzzy subsets. This is
accomplished by obtaining various relationships between the set of
all fuzzy ideals of a $\Gamma-$semiring and the set of all fuzzy
ideals of its left operator semiring such as lattice isomorphism
between the sets of fuzzy ideals of a $\Gamma-$semiring and its
operator semirings.
\end{abstract}

{\bf Keywards:} $\Gamma$-semiring, operator semiring, fuzzy
left(right) ideal, fuzzy k-ideal, $\Gamma$-semifield.

\section{Introduction}
The notion of $\Gamma-$semiring was introduced by M.M.K
Rao\cite{re:Rao}. This generalizes not only the notions of semiring
and $\Gamma-$ring but also the notion of ternary semiring. This also
provides an algebraic home to the nonpositive cones of the totally
ordered rings (It may be recalled here that the nonnegative cones of
the totally ordered rings form semirings but the nonpositive cones
do not as the induced multiplication is no longer closed).
$\Gamma-$semiring theory has been enriched with the help of operator
semirings of a $\Gamma-$semiring by Dutta and Sardar\cite{re:Dutta}.
To make operator semirings effective in the study of
$\Gamma-$semirings Dutta et al\cite{re:Dutta} established
correspondence between the ideals of a $\Gamma-$semiring S and the
ideals of the operator semirings of S. As it was done for
$\Gamma-$rings in \cite{re:Chanda} we also establish here various
relationships between the fuzzy ideals of a $\Gamma-$semiring S and
the fuzzy ideals of the operator semirings of S. We obtain a lattice
isomorphism between the sets of fuzzy ideals of a $\Gamma-$semiring
S and its left operator semiring L. This has been used to give a new
proof of the lattice isomorphism between the sets of ideals of S and
L which was originally due to Dutta and Sardar\cite{re:Dutta}. This
is also used to obtain a characterization of a $\Gamma-$semifield
and relationship between the fuzzy subsets of $\Gamma-$semiring S
and matrix $\Gamma_{n}-$semiring $S_{n}$.

\section{Preliminaries}
  We recall the following definitions and results for their use in
the sequel.

\begin{Definition} \textnormal{\cite{re:Zadeh}} Let S be a non empty set. A mapping
$\mu: S \rightarrow [0,1]$ is called a fuzzy subset of
S.\end{Definition}

\begin{Definition} \textnormal{\cite{re:Rao}}Let S and $\Gamma$ be two additive
commutative semigroups. Then S is called a $\Gamma-$semiring if
there exists a mapping $\\ S \times \Gamma \times S \rightarrow S $
(images to be denoted by $a \alpha b $ for $a, b \in S $ and $\alpha
\in \Gamma) $ satisfying the following conditions:\\
(i) $(a+b) \alpha c =a \alpha c+ b \alpha c$,\\
(ii) $ a \alpha (b+c)= a \alpha b+a \alpha c $,\\
(iii) $a (\alpha+\beta)b=a \alpha b +a \beta b$,\\
(iv) $ a \alpha (b \beta c)=(a \alpha b) \beta c $
for all $a,b,c \in S $ and for all $\alpha, \beta \in \Gamma $.\\
\\Further, if in a $\Gamma-$semiring, $(S,+)$ and $(\Gamma,+)$ are both
monoids and\\ (i) $ 0_{S} \alpha x=0_{S} =x \alpha 0_{S}$\\(ii)$ x
0_{\Gamma} y=0_{S}=y 0_{\Gamma} x$ for all $x,y \in S$ and for all
$\alpha \in \Gamma $ then we say that S is a $\Gamma-$semiring with
zero.\end{Definition}
 Throughout this paper we consider $\Gamma-$semiring with zero. For
simplification we write 0 instead of $0_{S}$.

\begin{Example} \textnormal{\cite{re:Dutta}}Let S be the additive commutative
semigroup of all $m \times n$ matrices over the set of all
non-negative integers and $\Gamma$ be the additive commutative
semigroup of all $n \times m$ matrices over the same set. Then S
forms a $\Gamma-$semiring if $x \alpha y$ denotes the usual matrix
multiplication of $x,\alpha,y$ where $x,y \in S$ and $\alpha \in
\Gamma$.\end{Example}

\begin{Definition} \textnormal{\cite{re:Dutta}}Let S be a $\Gamma-$semiring and F be the
free additive commutative semigroup generated by $ S \times \Gamma
$. Then the relation $\rho $ on F, defined by $
\displaystyle{\sum_{i=1}^{m}}(x_{i},\alpha_{i}) \rho
\displaystyle{\sum_{j=1}^{n}}(y_{j}, \beta_{j}) $ if and only if
$\displaystyle{\sum_{i=1}^{m}}x_{i} \alpha_{i} a=
\displaystyle{\sum_{j=1}^{n}}y_{j} \beta_{j} a $ \\for all $a \in S
(m,n \in Z^{+}) $, is a congruence on F. The congruence class
containing $\displaystyle{\sum_{i=1}^{m}}(x_{i},\alpha_{i})$ is
defined by $\displaystyle{\sum_{i=1}^{m}}[x_{i},\alpha_{i}]$. Then
$F/ \rho$, the set of congruence classes, is an additive commutative
semigroup and this also forms a semiring with the multiplication
defined by $$ \displaystyle{(\sum^{m}_{i=1}[x_{i},\alpha_{i}])
(\sum^{n}_{j=1}[y_{j},\beta_{j}])=\sum_{i,j}[x_{i}\alpha_{i}y_{j},\beta_{j}]}$$
\\This semiring is denoted by L and called the left operator
semiring of the $\Gamma-$semiring S.\end{Definition}

Dually the right operator semiring R of the
$\Gamma-$semiring S is defined.\\

\begin{Definition} \textnormal{\cite{re:Dutta}}Let S be a $\Gamma-$semiring and L be the
left operator semiring and R be the right one. If there exists an
element $\displaystyle{\sum^{m}_{i=1}[e_{i},\delta_{i}]} \in L
(\displaystyle{\sum^{n}_{j=1}}[\gamma_{j},f_{j}] \in R)$ such that
$\displaystyle{\sum^{m}_{i=1}e_{i}\delta_{i}a=a}
~~\displaystyle{(\sum^{n}_{j=1}a \gamma_{j}f_{j}=a)}$ for all $a \in
S $ then S is said to have the left unity
$\displaystyle{\sum^{m}_{i=1}[e_{i},\delta_{i}]}$ (resp. the right
unity
$\displaystyle{\sum^{n}_{j=1}[\gamma_{j},f_{j}]})$.\end{Definition}

\begin{Definition} \textnormal{\cite{re:Dutta}}Let S be a $\Gamma-$semiring and L be the
left operator semiring and R be the right one. For $ P\subseteq
L~(\subseteq R)$, $P^{+}:=\{a \in S: [a,\Gamma] \subseteq P\}$
(respectively $P^{*}:=\{a \in S:[\Gamma,a] \subseteq P\}$). For $Q
\subseteq S$, \\$Q^{+'}:=\{ \displaystyle{\sum_{i=1}^{m}} [x_{i},
\alpha_{i}] \in L :~(\displaystyle{\sum_{i=1}^{m}} ([x_{i},
\alpha_{i}])S\subseteq Q \}$ where ($\displaystyle{\sum_{i=1}^{m}}
[x_{i}, \alpha_{i}])S$ denotes the set of all finite sums
$\displaystyle{\sum_{i,k}}x_{i}\alpha_{i}s_{k},~~s_{k} \in S$ and\\
$Q^{*'}:=\{ \displaystyle{\sum_{i=1}^{m}} [\alpha_{i},x_{i}] \in R
:~(\displaystyle{\sum_{i=1}^{m}} S([ \alpha_{i},x_{i}])\subseteq Q
\}$ where $S(\displaystyle{\sum_{i=1}^{m}} [x_{i}, \alpha_{i}])$
denotes the set of all finite sums
$\displaystyle{\sum_{i,k}}s_{k}\alpha_{i}x_{i},~~s_{k} \in S$.
\end{Definition}
 Throughout this paper unless otherwise mentioned S denotes a $\Gamma$-semiring with left
unity and right unity and $ FLI(S) $, $ FRI(S) $ and $ FI(S) $
denote respectively the set of all fuzzy left ideals, the set of all
fuzzy right ideals and the set of all fuzzy ideals of the $ \Gamma
$-semiring S. Similar is the meaning of FLI(L), FRI(L), FI(L) where
L is the left operator semiring of the $\Gamma-$semiring S. Also
throughout we assume that $\mu(0)=1$ for a fuzzy left ideal (fuzzy
right ideal, fuzzy ideal) $\mu$ of a $\Gamma-$semiring S. Similarly
we assume that $\mu(0)=1$ for a fuzzy left ideal (fuzzy right ideal,
fuzzy ideal) $\mu$ of the left operator semring of a
$\Gamma-$semiring S.

\begin{Definition} \textnormal{\cite{re:Dutta}} A $\Gamma-$semiring S is said to be
zero-divisor free (ZDF) if $a \alpha b=0$ implies that either $a=0$
or $\alpha=0$ or $b=0$ for $a,b \in S,~\alpha \in
\Gamma$.\end{Definition}

\begin{Definition} \textnormal{\cite{re:Dutta}}A commutative $\Gamma-$semiring S is said to be a
$\Gamma-$semifield if for any $a (\neq 0)\in S$ and for any $\alpha
(\neq 0) \in \Gamma$ there exists $b \in S, \beta \in \Gamma$ such
that $a \alpha b \beta d =d$ for all $d \in S$.\end{Definition}

\begin{Example} \textnormal{\cite{re:Sardar}}Let $S=\{r\omega:r\in Q^{+}\cup \{0\}\}$
and $\Gamma=\{r\omega^{2}:r\in Q^{+}\cup
\{0\}\}$, where $Q^{+}$ is the set of all positive rational numbers.
Then S forms a $\Gamma-$semifield.\end{Example}

\begin{Definition} \textnormal{\cite{re:Biswas1}}Let $\mu$ be a non empty fuzzy subset of a
semiring S (i.e. $\mu(x) \neq 0 $ for some $x \in S
 $). Then $\mu $ is called a fuzzy left ideal [ fuzzy right ideal]
 of S if \\(i) $\mu(x+y) \geq min [\mu(x), \mu(y)]$ and
  \\~~~~~~~~(ii) $\mu(xy) \geq \mu(y) $ [resp. $\mu(xy) \geq
  \mu(x)]$ for all $ x,y \in S$.\end{Definition}

\begin{Definition} \textnormal{\cite{re:Goswami}}Let $\mu$ be a non empty fuzzy subset of a
$\Gamma-$semiring S (i.e. $\mu(x) \neq 0 $ for some $x \in S
 $). Then $\mu $ is called a fuzzy left ideal [ fuzzy right ideal]
 of S if \\(i) $\mu(x+y) \geq min [\mu(x), \mu(y)]$ and
  \\~~~~~~~~(ii) $\mu(x \gamma y) \geq \mu(y) $ [resp. $\mu(x \gamma y) \geq
  \mu(x)]$ for all $ x,y \in S, \gamma \in \Gamma$.\end{Definition}

\begin{Definition} \textnormal{\cite{re:Goswami2}}Let S be a $ \Gamma $-semiring and $
\mu_{1}, \mu_{2} \in FLI(S) ~[ FRI(S),~FI(S)] $. Then the sum $
\mu_{1}\oplus \mu_{2}$ of $ \mu_{1} $ and $\mu_{2} $ is defined as
follows:

$ (\mu_{1} \oplus \mu_{2})(x)= \displaystyle{\sup_{x=u+v}}[\min[
\mu_{1}(u), \mu_{2}(v) ]: u, v \in S]\\~~~~~~~~~~~~~~~~~~~~ = 0$ if
for any $ u,v \in S, u+v \neq x $.\end{Definition}

\textbf{ Note.} Since S contains $ 0 $, in the above definition the
case $ x \neq u+v $ for any $ u,v \in S $ does not arise.

\begin{Definition} \textnormal{\cite{re:Matrix}}Let S be a $ \Gamma
$-semiring and n be a positive integer. The sets of $n\times n$
matrices with entries from S and $n\times n$ matrices with entries
from $\Gamma$ are denoted by $S_{n}$ and $\Gamma_{n}$ respectively.
Let $A,B \in S_{n}$ and $\Delta \in \Gamma_{n}$. Then $A\Delta B \in
S_{n}$ and $A+B \in S_{n}$. Clearly, $S_{n}$ forms a
$\Gamma_{n}$-semiring with these operations. This is called the
matrix $\Gamma-$semiring over S or the matrix $\Gamma_{n}-$semiring
$S_{n}$ or simply the $\Gamma_{n}-$semiring $S_{n}$.\end{Definition}

  The right operator semiring of the matrix $\Gamma_{n}-$semiring
  $S_{n}$ is denoted by $[\Gamma_{n},S_{n}]$ and the left one by
  $[S_{n},\Gamma_{n}]$. If $x \in S$, the notation $xE_{ij}$ will be
  used to denote a matrix in $S_{n}$ with $x$ in the (i,j)-th entry
  and zeros elsewhere. $\alpha E_{ij} \in \Gamma_{n}$, where $\alpha \in
  \Gamma$, will have a similar meaning. If $P \subseteq S,~P_{n}$
  will denote the set of all $n\times n$ matrices with entries from
  P. Similar is the meaning of $\Delta_{n}$ where $\Delta \subseteq
  \Gamma$.

\begin{Proposition} \textnormal{\cite{re:Matrix}}Let S be a
$\Gamma-$semiring and R (resp. L) be its right (resp. left) operator
semiring. Let $R_{n}$ (resp. $L_{n}$) denote the semiring of all
$n\times n$ matrices over R (resp. L). Then
\\(i) the right operator semiring $[\Gamma_{n}, S_{n}]$ of the
$\Gamma_{n}-$semiring $S_{n}$ is isomorphic with $R_{n}$ via the
mapping $\displaystyle{\sum_{i=1}^{p}[[\gamma_{jk}^{i}],[x_{uv}^{i}]
]\mapsto
\sum_{i=1}^{p}(\sum_{t=1}^{n}[\gamma_{jt}^{i},x_{tv}^{i}])_{1 \leq
j,v\leq n}}$,
\\and (ii) the left operator semiring $[S_{n},\Gamma_{n}]$ of the
$\Gamma_{n}-$semiring $S_{n}$ is isomorphic with $L_{n}$ via the
mapping $\displaystyle{\sum_{i=1}^{p}[[x_{uv}^{i}],[\gamma_{jk}^{i}]
]\mapsto
\sum_{i=1}^{p}(\sum_{t=1}^{n}[x_{ut}^{i},\gamma_{tk}^{i}])_{1 \leq
u,k\leq n}}$.\label{Proposition:2.10} \end{Proposition}

  In view of the above proposition we henceforth identify
$[S_{n},\Gamma_{n}]$ and $R_{n}$; $[S_{n},\Gamma_{n}]$ and $L_{n}$.

\section{Main results.}
 Throughout this section S denotes a $\Gamma-$semiring and L denotes the left operator
 semiring of the $ \Gamma-$semiring S.

\begin{Definition} Let $ \mu$ be a fuzzy subset of L, we
define a fuzzy subset $\mu^{+}$ of S by $ \mu^{+}(x)= \displaystyle
{\inf_{\gamma \in \Gamma}} \mu([x, \gamma]) $ where $x \in S $.
\\If $\sigma$ is a fuzzy subset of S, we define a fuzzy
subset $\sigma^{+'}$ of L by $\\
\sigma^{+'}(\displaystyle{\sum_{i}}[x_{i}, \alpha_{i}])
=\displaystyle{\inf_{s \in S}} \sigma(\displaystyle{\sum_{i}} x_{i}
\alpha_{i}s) $ where $ \displaystyle{\sum_{i}} [x_{i}, \alpha_{i}]
\in L $. \label{Def:3.1}\end{Definition}

By routine verification we obtain the following lemma.

\begin{Lemma} If $\{ \mu_{i}: i \in I \}$ is a collection of
fuzzy subsets of L then $\displaystyle{ \bigcap_{i \in I}
\mu_{i}^{+}=(\bigcap_{i \in I}
\mu_{i})^{+}}$.\label{Lemma:3.3}\end{Lemma}

\begin{Proposition} Suppose $\sigma,\sigma_{1}, \sigma_{2} \in FI(S)$and $\mu \in
FI(L)$. Then
\\(i) $\sigma^{+'} \in FI(L)$. Moreover, if $\sigma$ is non constant
then $\sigma^{+'}$ is non constant.\label{Proposition:3.4(i)}
\\(ii) $(\sigma^{+'})^{+}=\sigma $,
\\(iii) $\sigma_{1} \neq \sigma_{2}$ implies that $\sigma_{1}^{+'} \neq
\sigma_{2}^{+'}$,
\\(iv) $(\sigma_{1} \oplus \sigma_{2})^{+'}=\sigma^{+'}_{1} \oplus
\sigma^{+'}_{2}$,
\\(v) $(\sigma_{1} \cap \sigma_{2})^{+'}=\sigma_{1}^{+'} \cap
\sigma_{2}^{+'}$,
\\(vi) $\sigma_{1} \subseteq \sigma_{2}$ implies that $\sigma_{1}^{+'} \subseteq
\sigma_{2}^{+'}$,
\\(vii) $\mu^{+} \in FI(S)$. Moreover, if $\mu$ is non constant then
$\mu^{+}$ is non constant.\label{Proposition:3.4vii)}
\\(viii) $(\mu^{+})^{+'} = \mu $,
\\(ix) $\mu_{1} \subseteq \mu_{2}$ implies that $\mu_{1}^{+} \subseteq
\mu_{2}^{+}$.\label{Proposition:3.4}\end{Proposition}

\textit{Proof.} (i) Let $\sigma \in FI(S)$. Then $\sigma(0)=1$. Let
$\gamma \in \Gamma$. Then \\$\sigma^{+'}([0, \gamma])=
\displaystyle{\inf_{s \in S} [\sigma(0 \gamma s)]
=\sigma(0)=1,~~~\gamma \in \Gamma}$. Thus we see that $\sigma^{+'}$
is non empty and $\sigma^{+'}(0)=1$ as for all $\gamma \in \Gamma,~
[0,\gamma]$ is the zero element of L.
\\Let $\displaystyle{\sum_{i} [x_{i}, \alpha_{i}],~~\sum_{j} [y_{j}, \beta_{j}]} \in L $.
\\Then $\sigma^{+'}\displaystyle{(\sum_{i}[x_{i}, \alpha_{i}]+ \sum_{j}[y_{j},
\beta_{j}])} =\displaystyle{\inf_{s \in S}[\sigma(\sum_{i} x_{i}
\alpha_{i}s +\sum_{j}y_{j} \beta_{j}s)]}
\\~~~~~~~~~~~\geq \displaystyle{\inf_{s \in S}[\min[\sigma(\sum_{i} x_{i} \alpha_{i}s),
 \sigma(\sum_{j} y_{j} \beta_{j}s)]]}
\\~~~~~~~~~~~=\min[\displaystyle{\inf_{s \in S} [\sigma(\sum_{i} x_{i}
\alpha_{i}s)], \inf_{s \in S}[\sigma(\sum_{j} y_{j} \beta_{j}
s)]}]\\~~~~~~~~~~~=\min[ \sigma^{+'}(\displaystyle{\sum_{i}[x_{i},
\alpha_{i}])}, \sigma^{+'}(\displaystyle{\sum_{j}[y_{j},
\beta_{j}])}]$.
\\Again $\sigma^{+'}(\displaystyle{\sum_{i}[ x_{i}, \alpha_{i}] \sum_{j}[ y_{j}, \beta_{j}]})
=\sigma^{+'}\displaystyle{(\sum_{i,j}[ x_{i} \alpha_{i} y_{j},
\beta_{j}])} =\displaystyle{\inf_{s \in S}}
\sigma(\displaystyle{\sum_{i,j}} x_{i} \alpha_{i} y_{j} \beta_{j} s)
 \\~~~~~~~~~~~~\geq \displaystyle{
 [\min[ \sigma(\sum_{i} x_{i} \alpha_{i}y_{1}),
\sigma(\sum_{i}~~ x_{i} \alpha_{i}y_{2}), \sigma(\sum_{i} x_{i}
\alpha_{i}y_{3}),..........]]}\\~~~~~~~~~~~~\geq
\displaystyle{\inf_{s \in S}[ \sigma(\sum_{i}(x_{i} \alpha_{i}
s))]}=\sigma^{+'}(\displaystyle{\sum_{i}}[x_{i}, \alpha_{i}])$.
\\Similarly we can show that $\displaystyle{\sigma^{+'}( \sum_{i}[ x_{i}, \alpha_{i}]\sum_{j}[ y_{j},
\beta_{j}]) \geq \sigma^{+'}(\sum_{j}[y_{j},\beta_{j}])}$.
\\Consequently, $\sigma^{+'}\in FI(L)$.
\\Further, let $\sigma$ be a non constant fuzzy ideal of S.
\\If $\sigma^{+'}$ is constant then
$\sigma^{+'}(\displaystyle{\sum_{i}[x_{i},\alpha_{i}]})=1$ for all
$\displaystyle{\sum_{i}[x_{i},\alpha_{i}]} \in L$. Then
$\sigma(x)=(\sigma^{+'})^{+}(x)=\displaystyle{\inf_{\gamma \in
\Gamma}}\sigma^{+'}([x,\gamma])=1$ for all $x \in S$ where $\sigma$ is constant. This contradicts that $\sigma$ is non constant . Consequently, $\sigma^{+'}$ is non constant.
\\(ii) Let $x \in S $. Then $\\((\sigma^{+'})^{+})(x)=
\displaystyle{\inf_{\gamma \in \Gamma} [\sigma^{+'}([x,
\gamma])]=\inf_{\gamma\in \Gamma}[ \inf_{s \in S }[\sigma(x \gamma
s)]]} \geq \sigma(x)$. Hence $\sigma \subseteq (\sigma^{+'})^{+}$.
\\Let $\displaystyle{\sum_{i}} [\gamma_{i}, f_{i}]$ be the right unity of S.
Then $\displaystyle{\sum_{i}} x \gamma_{i} f_{i}=x $ for all $x \in
S$.
\\Now $\sigma(x)=\sigma(\displaystyle{\sum_{i}} x \gamma_{i} f_{i})
\geq \min [\sigma(x \gamma_{1} f_{1}),\sigma(x \gamma_{2}
f_{2}),.........]\\~~~~~~~~~~~~~~\geq \displaystyle{\inf_{\gamma \in
\Gamma}[ \inf_{s \in S} [\sigma(x \gamma s)] ]}=
(\sigma^{+'})^{+}(x)$.
\\Therefore $(\sigma^{+'})^{+} \subseteq \sigma $ and hence $ (\sigma^{+'})^{+}=\sigma
$.
\\(iii) Let $\sigma_{1} \neq \sigma_{2}$. If possible, let
$\sigma_{1}^{+'}=\sigma_{2}^{+'}$. Then $(\sigma_{1}^{+'})^{+} =
(\sigma_{2}^{+'})^{+}$. i.e., $\sigma_{1} =\sigma_{2}$, (by (ii)) --
which contradicts our assumption. Hence $\sigma_{1}^{+'} \neq
\sigma_{2}^{+'}$
\\(iv) Let $\displaystyle{ \sum_{i}[a_{i}, \alpha_{i}]} \in L $. Then
\\ $((\sigma_{1} \oplus \sigma_{2})^{+'})(\displaystyle{\sum_{i}}[a_{i},
\alpha_{i}])$=$\displaystyle{\inf_{s \in S}(\sigma_{1} \oplus
\sigma_{2})(\sum_{i} a_{i} \alpha_{i}s)}\\= \displaystyle{\inf_{s
\in S}[\sup[\min[ \sigma_{1}(\sum_{k} x_{k} \delta_{k}s),
\sigma_{2}(\sum_{j} y_{j} \beta_{j}s)]: \sum_{i}a_{i}
\alpha_{i}s=\sum_{k} x_{k} \delta_{k}s+\sum_{j} y_{j} \beta_{j}s]]}$
\\ $~~~~~~=\displaystyle{\sup[\min[\inf_{s \in S} \sigma_{1}(\sum_{k} x_{k} \delta_{k}s),
\inf_{s \in S} \sigma_{2}(\sum_{j}y_{j} \beta_{j}s)]]}$
\\ $~~~~~~=\displaystyle{\sup[\min[\sigma^{+'}_{1}(\sum_{k}[x_{k}, \delta_{k}]),
 \sigma^{+'}_{2}(\sum_{j}[y_{j}, \beta_{j}])]]}\\~~~~~~=(\sigma^{+'}_{1} \oplus \sigma^{+'}_{2})
 (\displaystyle{\sum_{i}}[a_{i}, \alpha_{i}])$
\\Thus $(\sigma_{1} \oplus \sigma_{2})^{+'}=\sigma^{+'}_{1} \oplus
\sigma^{+'}_{2}$.
\\(v) $\displaystyle{(\sigma_{1} \cap \sigma_{2})^{+'}(\sum_{i}[a_{i},
\alpha_{i}])= \inf_{s \in S}[(\sigma_{1} \cap \sigma_{2})(\sum_{i}
a_{i}\alpha_{i}s)]}\\~~~~~~~~~~~~~~~~=\displaystyle{\inf_{s \in
S}[\min[\sigma_{1}(\sum_{i}a_{i}\alpha_{i}s),\sigma_{2}(\sum_{i}a_{i}\alpha_{i}s)]]}
\\~~~~~~~~~~~~~~~~=\displaystyle{\min[\inf_{s \in S}\sigma_{1}(\sum_{i}a_{i}\alpha_{i}s),
\inf_{s \in S}\sigma_{2}(\sum_{i}a_{i}\alpha_{i}s)]}
\\~~~~~~~~~~~~~~~~=\displaystyle{\min[\sigma_{1}^{+'}(\sum_{i}[a_{i},\alpha_{i}]),
\sigma_{2}^{+'}(\sum_{i}[a_{i},\alpha_{i}])]}
\\~~~~~~~~~~~~~~~~=(\sigma_{1}^{+'} \cap
\sigma_{2}^{+'})(\displaystyle{\sum_{i}}[a_{i}, \alpha_{i}]) $
\\Hence $(\sigma_{1} \cap \sigma_{2})^{+'}=\sigma_{1}^{+'} \cap
\sigma_{2}^{+'}$.
\\(vi) Let $ \sigma_{1}, \sigma_{2} \in FI(S)$ be such that
$\sigma_{1} \subseteq \sigma_{2}$. Then
\\ $\displaystyle{\sigma_{1}^{+'}(\sum_{i}[x_{i}, \alpha_{i}])}= \displaystyle{\inf_{s \in S}[\sigma_{1}(\sum_{i} x_{i}
\alpha_{i}s)]} \leq \displaystyle{\inf_{s \in S}[\sigma_{2}(\sum_{i}
x_{i} \alpha_{i}s)]}
\\~~~~~~~~~~~~~~~~~~~~=\sigma_{2}^{+'}(\displaystyle{\sum_{i}}[x_{i}, \alpha_{i}])$ for all
$\displaystyle{\sum_{i}}[x_{i}, \alpha_{i}] \in L $.
\\Thus $\sigma_{1}^{+'} \subseteq \sigma_{2}^{+'}$.
\\(vii) Let $\mu \in FI( L)$. Then $\mu(0)=1$.
\\Now $\mu^{+}(0)= \displaystyle{\inf_{\gamma \in \Gamma} [\mu([0,
\gamma])]} =1$ [Since for all $\gamma \in \Gamma ,~[0, \gamma] $ is
the zero element of L].
\\This also shows that $\mu^{+}$ is non empty.
\\Let $x, y \in S $ and $\alpha \in \Gamma$.
Then\\ $\mu^{+}(x+y)=\displaystyle{\inf_{\gamma \in
\Gamma}[\mu([x+y, \gamma])]} =\inf[\mu([x, \gamma]+[y, \gamma])]\geq
\inf[\min[\mu([x, \gamma]),
\mu([y,\gamma])]]\\=\min[\displaystyle{\inf_{\gamma \in \Gamma}
[\mu([x, \gamma])], \inf_{\gamma \in \Gamma} [\mu([y,
\gamma])]}=\min[\mu^{+}(x), \mu^{+}(y)]$
\\Therefore $\mu^{+}(x+y) \geq \min[\mu^{+}(x), \mu^{+}(y)]$.
\\Again $\mu^{+}(x \alpha y)=\displaystyle{\inf_{\gamma \in \Gamma}[\mu([x \alpha y, \gamma])]}
= \displaystyle{\inf_{\gamma \in \Gamma}[\mu([x, \alpha][y,
\gamma])]}\geq \displaystyle{\inf_{\gamma \in \Gamma} \mu[y,
\gamma]}=\mu^{+}(y)$ and  $\mu^{+}(x \alpha
y)=\displaystyle{\inf_{\gamma \in \Gamma}[\mu([x \alpha
y,\gamma])]}=\displaystyle{\inf_{\gamma \in \Gamma}[\mu([x,
\alpha][y, \gamma])]}\geq \mu([x, \alpha]) \geq
\displaystyle{\inf_{\delta \in \Gamma}[\mu([x,
\delta])]}\\~~~~~~~~~~~~~~~~~~=\mu^{+}(x)$.
\\Consequently $\mu^{+}\in FI(S)$.
\\Suppose $\mu$ is a non constant fuzzy ideal of L.
If $\mu^{+}$ is constant then $\mu^{+}(x)=1$ for all $x \in S$ as $\mu^{+}(0_{S})=1$. Suppose\\
$\displaystyle{\sum_{i}}[x_{i},\alpha_{i}] \in L$. Then
$\mu(\displaystyle{\sum_{i}}[x_{i},\alpha_{i}])=((\mu^{+})^{+'})(\displaystyle{\sum_{i}}[x_{i},\alpha_{i}])
=\displaystyle{\inf_{s \in S}\mu^{+}(\sum_{i}x_{i}\alpha_{i}
s)}\\=1$ for all $\displaystyle{\sum_{i}}[x_{i},\alpha_{i}] \in L$.
This implies that $\mu$ is constant, which contradicts our
assumption. This completes the proof.
\\(viii) Let $\mu \in FI(L)$. Then for $\displaystyle{\sum_{i}}[x_{i}, \alpha_{i}]\in
L$,
\\$ ((\mu^{+})^{+'})(\displaystyle{\sum_{i}}[x_{i}, \alpha_{i}])
=\displaystyle{\inf_{s \in S}[\mu^{+}(\sum_{i} x_{i} \alpha_{i} s)]}
=\displaystyle{\inf_{s \in S}[\inf_{\gamma \in \Gamma}[\mu(\sum_{i}
[x_{i}
 \alpha_{i} s, \gamma])]]}
\\~~~~~~~~~~~~~~~~~~~~~~~~~~~=\displaystyle{\inf_{s \in S}[\inf_{\gamma \in \Gamma}[\mu(\sum_{i}
[ x_{i}, \alpha_{i}][ s, \gamma])]]} \geq\mu(\sum_{i} [x_{i},
\alpha_{i}])$.
\\Hence $\mu \subseteq (\mu^{+})^{+'}$.
\\Let $\displaystyle{\sum_{i}}[e_{i}, \delta_{i}]$ be the left unity of S.Then
\\ $ \displaystyle{\mu(\sum_{j}[x_{j}, \alpha_{j}])= \mu(\sum_{j}[x_{j}, \alpha_{j}]~\sum_{i}[e_{i},
\delta_{i}])}
\\~~~~~~~~~~~~~\geq \displaystyle{\min[\mu(\sum_{j}[x_{j}, \alpha_{j}][e_{1}, \delta_{1}]),
\mu(\sum_{j}[x_{j}, \alpha_{j}][e_{2},
\delta_{2}]),............]}\\~~~~~~~~~~~\geq \displaystyle{\inf_{ s
\in S}[\inf_{\gamma \in \Gamma}[\mu(\sum_{j}[x_{j}, \alpha_{j}][s,
\gamma]) ]]}=(\mu^{+})^{+'}(\displaystyle{\sum_{j}}[x_{j},
\alpha_{j}])$.
\\Therefore $(\mu^{+})^{+'} \subseteq \mu $ and so $(\mu^{+})^{+'} = \mu $.
\\(ix) Proof is similar to that of (vi).

\textbf{Note.} All the results of the above proposition also hold
for FRI(S) and FRI(R).
\\ In view of the above Proposition and the fact that the set of fuzzy ideals of a
$\Gamma-$semiring S and that of its operator semirings form a
lattice under the operations $\oplus$ and $\cap$ we obtain the
following result.

\begin{Theorem} The lattices of all fuzzy ideals [ fuzzy right
ideals ] of S and L are isomorphic via the inclusion preserving
bijection $\sigma \mapsto \sigma^{+'}$ where $\sigma \in FI(S) $
[resp. FRI(S)] and $\sigma^{+'} \in FI(L) $ [resp.
FRI(L)].\label{Th:3.8}\end{Theorem}

\begin{Corollary} FLI(S) [resp. FRI(S), FI(S)]] is a complete
lattice.\label{Cor:3.10}\end{Corollary}

\textit{Proof.} The corollary follows from the above theorem and the
fact that FLI(L) [resp. FRI(L), FI(L)] is a complete
lattice\cite{re:Biswas}.

\begin{Lemma} Let I be an ideal (left ideal, right ideal) of a $\Gamma-$semiring S and
$\lambda_{I}$ be the characteristic function of I. Then
$(\lambda_{I})^{+'}=\lambda_{I^{+'}}$. Moreover, $I^{+'}$ is an ideal of L. \label{Lemma:3.11}\end{Lemma}

\textit{Proof.} Let
$\displaystyle{\sum_{i}}[x_{i},\alpha_{i}] \in L$. Then either
$\displaystyle{\sum_{i}}[x_{i},\alpha_{i}] \in I^{+'}$ or
$\displaystyle{\sum_{i}}[x_{i},\alpha_{i}] \notin I^{+'}$. If
$\displaystyle{\sum_{i}}[x_{i},\alpha_{i}] \in I^{+'}$ then
$\lambda_{I^{+'}}(\displaystyle{\sum_{i}}[x_{i},\alpha_{i}])=1$ and
\\$\displaystyle{\sum_{i}}x_{i} \alpha_{i}s \in I $ for all $s \in
S$. Hence
$(\lambda_{I})^{+'}(\displaystyle{\sum_{i}}[x_{i},\alpha_{i}]
)=\displaystyle{\inf_{s \in S}}
\lambda_{I}(\displaystyle{\sum_{i}}x_{i} \alpha_{i}s)=1$.
\\Again, if $\displaystyle{\sum_{i}}[x_{i},\alpha_{i}] \notin I^{+'}$
then $\lambda_{I^{+'}}(\displaystyle{\sum_{i}}[x_{i},\alpha_{i}])=0$
and $\displaystyle{\sum_{i}}x_{i} \alpha_{i}s \notin I $ for some $s
\in S$. Hence
$(\lambda_{I})^{+'}(\displaystyle{\sum_{i}}[x_{i},\alpha_{i}])
=\displaystyle{\inf_{s \in S}}
\lambda_{I}(\displaystyle{\sum_{i}}x_{i} \alpha_{i}s)=0$. Thus we
obtain $(\lambda_{I})^{+'}=\lambda_{I^{+'}}$.
Now by Proposition \ref{Proposition:3.4}(i), $(\lambda_{I})^{+'}$ is a fuzzy ideal of the left operator semiring L. Hence $\lambda_{I^{+'}}$ is a fuzzy ideal of L. Consequently, $I^{+'}$ is an ideal of L\cite{re:Biswas1}

\begin{Lemma} Let I be an ideal (left ideal, right ideal) of the left
operator semiring L of a $\Gamma-$semiring S and $\lambda_{I}$ be
the characteristic finction of I. Then
$(\lambda_{I})^{+}=\lambda_{I^{+}}$. Moreover, $I^{+}$ is an ideal of S. \label{Lemma:3.12}\end{Lemma}

\textit{Proof.}  Let $x \in S$. Then
either $x \in I^{+} $ or $x \notin I^{+}$. If $x \in I^{+}$ then
$\lambda_{I^{+}}(x)=1$. Again $x \in I^{+}$ implies that $[x,
\gamma] \in I$ for all $\gamma \in \Gamma$. Hence
$(\lambda_{I})^{+}(x)=\displaystyle{\inf_{\gamma \in \Gamma}
\lambda_{I}([x,\gamma])}=1$.
\\Again, if $x \notin I^{+}$ then $\lambda_{I^{+}}(x)=0$ and
$[x, \gamma] \notin I$ for some $\gamma \in \Gamma$. Hence
$(\lambda_{I})^{+}(x)=\displaystyle{\inf_{\gamma \in \Gamma}
\lambda_{I}([x,\gamma])}=0$. Consequently,
$(\lambda_{I})^{+}=\lambda_{I^{+}}$. Now by Proposition \ref{Proposition:3.4}(ii), $(\lambda_{I})^{+}$ is a fuzzy ideal of S. Hence $I^{+}$ is an ideal of S\cite{re:Goswami}

\textbf{Remark.} The last parts of Lemma \ref{Lemma:3.11} and Lemma \ref{Lemma:3.12} are originally due to Dutta and Sardar. They are established here via fuzzy subsets.

The following theorem is also due to Dutta and Sardar\cite{re:Dutta}. We
give an alternative proof of it by using the lattice isomorphism of fuzzy
ideals obtained in Theorem \ref{Th:3.8} and by using Lemma \ref{Lemma:3.11}
and Lemma \ref{Lemma:3.12}.

\begin{Theorem} The lattices of all ideals [ right
ideals ] of S and L are isomorphic via the mapping $I \mapsto
I^{+'}$ where I denotes an ideal ( right ideal) of
S.\label{Th:3.15}\end{Theorem}

\textit{Proof.} That $I \mapsto I^{+'}$ is a mapping follows from
Lemma \ref{Lemma:3.11}. Let $I_{1}$ and $I_{2}$ be two
ideals of S such that $I_{1} \neq I_{2}$. Then $\lambda_{I_{1}}$ and
$\lambda_{I_{2}}$ are fuzzy ideals of S where $\lambda_{I_{1}}$ and
$\lambda_{I_{2}}$ are characteristic functions of $I_{1}$ and
$I_{2}$ respectively\cite{re:Goswami}. Evidently, $\lambda_{I_{1}}
\neq \lambda_{I_{2}}$. Then by Theorem \ref{Th:3.8},
$\lambda_{I_{1}}^{+'} \neq \lambda_{I_{2}}^{+'}$. Hence by Lemma
\ref{Lemma:3.11}, $\lambda_{I_{1}^{+'}} \neq \lambda_{I_{2}^{+'}}$
whence $I_{1}^{+'} \neq I_{2}^{+'}$ Consequently, the mapping $I
\mapsto I^{+'}$ is one-one.
\\Next let J be an ideal of L. Then $\lambda_{J}$ is a fuzzy ideal
of L\cite{re:Biswas}. By Proposition \ref{Proposition:3.4vii)}(vii),
$(\lambda_{J})^{+}$ is a fuzzy ideal of S and by Theorem
\ref{Th:3.8} we obtain $((\lambda_{J})^{+})^{+'}=\lambda_{J}$ . Now
by successive use of Lemma \ref{Lemma:3.12} and Lemma \ref{Lemma:3.11} we obtain
$\lambda_{({J}^{+})^{+'}}=\lambda_{J}$ and consequently,
$(J^{+})^{+'}=J$. Hence the mapping is onto.
\\Let $I_{1},~I_{2}$ be two ideals of S such that $I_{1} \subseteq
I_{2}$. Then $\lambda_{I_{1}} \subseteq \lambda_{I_{2}}$. Hence by
Theorem \ref{Th:3.8}, $\lambda_{I_{1}}^{+'} \subseteq
\lambda_{I_{2}}^{+'}$ whence by Lemma \ref{Lemma:3.11},
$\lambda_{I_{1}^{+'}} \subseteq \lambda_{I_{2}^{+'}}$. Consequently
$I_{1}^{+'} \subseteq I_{2}^{+'}$. Thus the mapping is inclusion
preserving. This completes the proof.

\begin{Theorem} A commutative semiring M is a semifield if
and only if for every non constant fuzzy ideal $\mu$ of M,
$\mu(x)=\mu(y) < \mu(0)$ for all $x,y \in M\setminus \{
0\}$.\label{Th:3.17}\end{Theorem}

\textit{Proof.} Let M be a semifield and let $x,y \in M\setminus \{0
\}$ and $\mu$ a non constant fuzzy ideal of M. Then
$\mu(y)=\mu(yxx^{-1}) \geq \mu(xx^{-1}) \geq \mu(x)$ where $x^{-1}$
is the inverse of $ x $. Similarly $\mu(x) \geq \mu(y)$. Therefore
$\mu(x)=\mu(y)$. It is known that $\mu(0) \geq \mu(x)$ for all $x
\in M$.
  Now we claim that $\mu(0) > \mu(x)$
for all $x \in M\setminus \{0 \}$. Otherwise, if $\mu(0) = \mu(x)$
for some $x \in M\setminus \{0 \}$ then by what we have obtained,
$\mu(0) = \mu(x)$ for all $ x \in M\setminus \{ 0 \}$ --a
contradiction to the supposition that $\mu$ is non constant.
\\Thus $\mu(x)=\mu(y) < \mu(0)$ for all $x,y \in M\setminus \{ 0\}$.

Conversely, let M be a commutative semiring and for every non
constant fuzzy ideal $\mu$ of M, $\mu(x)=\mu(y) < \mu(0)$ for all
$x,y \in M\setminus \{ 0\}$. Now let I be a non zero ideal of M. If
possible, suppose $I\neq M$. Then there exists an element $x \in M
\setminus I$.
\\Let $\lambda_{I}$ be the characteristic function of I. Then $\lambda_{I}(x)=0 \neq
1=\lambda_{I}(0)$. This implies that $\lambda_{I}$ is a non constant
fuzzy ideal of M. Suppose $y (\neq 0) \in I$. Then by hypothesis
$\lambda_{I}(x)=\lambda_{I}(y)$ but $\lambda_{I}(x)=0 $ and $
\lambda_{I}(y)=1$. Thus we get a contradiction. Hence $I=M$. Thus M
has no non-zero proper ideals. Consequently, M is a semifield\cite{re:Golan}.

\begin{Theorem} A ZDF commutative $\Gamma-$semiring S is a
$\Gamma-$semifield if and only if for every non constant fuzzy ideal
$\mu$ of S, $\mu(x)=\mu(y) < \mu(0)$ for all $x,y \in S\setminus \{
0\}$.\label{Th:3.18}\end{Theorem}

\textit{Proof.} Let S be a ZDF commutative $\Gamma-$semifield and
$\mu$ a non constant fuzzy ideal of S. Let $x,y \in S\setminus \{0
\}$ and $\alpha(\neq 0) \in \Gamma$. Then there exists $z \in S,
\beta \in \Gamma$ such that $x \alpha z \beta s=s$ for all $s \in
S$. In particular, $x \alpha z \beta y=y$. Then $\mu(y)=\mu(x \alpha
z \beta y)\geq \mu(x)$. Similarly $\mu(x) \geq \mu(y)$. Therefore
$\mu(x)=\mu(y)$.
\\  Now we claim that $\mu(0) > \mu(x)$
for all $x \in S\setminus \{0\}$. Otherwise, if $\mu(0) = \mu(x)$
for some $x \in S\setminus \{0\}$ then $\mu(0) = \mu(x)$ for all $ x
\in S\setminus \{ 0\}$ --a contradiction that $\mu$ is non constant.
\\Thus $\mu(x)=\mu(y) < \mu(0)$ for all $x,y \in S\setminus \{ 0\}$.
\\Converse follows by applying the similar argument as applied in
the converse part of Theorem \ref{Th:3.17} and by using Theorem 9.7
of \cite{re:Dutta}.

Using the above theorem and Proposition \ref{Proposition:3.4} (hence
Theorem \ref{Th:3.8}) we give a new proof of the following result of Dutta and Sardar\cite{re:Dutta}.

\begin{Theorem} Let S is a ZDF commutative $\Gamma-$semiring. Then S is a $\Gamma-$semifield if and only if
 its left operator semiring L is a semifield.\end{Theorem}

\textit{Proof.} Let S is a $\Gamma-$semifield. Let $\mu$ be a non
constant fuzzy ideal of L. Then by Proposition
\ref{Proposition:3.4}(vii), $\mu^{+}$ is a non constant fuzzy ideal
of S. Hence by Theorem \ref{Th:3.18}, $\mu^{+}(x)=\mu^{+}(y) <
\mu^{+}(0_{S})$ for all $x,y
\in S\setminus \{ 0_{S}\}$. Let\\
$\displaystyle{\sum_{i}}[x_{i},\alpha_{i}],
\displaystyle{\sum_{j}}[x_{j},\beta_{j}]  \in L\setminus \{
0_{L}\}$. Then
$\mu(\displaystyle{\sum_{i}}[x_{i},\alpha_{i}])=((\mu^{+})^{+'})(\displaystyle{\sum_{i}}[x_{i},\alpha_{i}])
\\=\displaystyle{\inf_{s \in S}\mu^{+}(\sum_{i}x_{i}\alpha_{i}}
s)=\displaystyle{\inf_{s \in
S}\mu^{+}(\sum_{j}x_{j}\beta_{j}s)}=((\mu^{+})^{+'})(\displaystyle{\sum_{j}}[x_{j},\beta_{j}])
< \mu(0_{L})$.\\Hence by Theorem \ref{Th:3.17}, L is a semifield.
\\Conversely, let L be a semifield and let $\mu$ be a non
constant fuzzy ideal of S. Then by Proposition
\ref{Proposition:3.4}(i), $\mu^{+'}$ is a non constant fuzzy ideal
of L. Hence by Theorem
\ref{Th:3.17}, $\mu^{+'}(\displaystyle{\sum_{i}[x_{i},\alpha_{i}]})=
\mu^{+'}(\displaystyle{\sum_{j}[y_{j},\alpha_{j}]})
<\mu^{+'}(0_{L})$ for all
\\$\displaystyle{\sum_{i}[x_{i},\alpha_{i}]},\displaystyle{\sum_{j}[y_{j},\alpha_{j}]}
\in L\setminus\{0_{L}\}$. Now let $x,y \in S\setminus\{0_{S}\}$.
Then\\$\mu(x)=(\mu^{+'})^{+}(x)=\displaystyle{\inf_{\gamma \in
\Gamma}}\mu^{+'}([x,\gamma])=\displaystyle{\inf_{\gamma \in
\Gamma}}\mu^{+'}([y,\gamma])=\mu(y)<\mu(0_{S})$. Hence S is a
$\Gamma-$semifield.

The following result is also a consequence of Theorem \ref{Th:3.8}.

\begin{Theorem} Suppose S is a $\Gamma-$semiring with unities and n is a positive
integer. Then there exists an inclusion preserving bijection between
the set of all fuzzy ideals of S and the set of all fuzzy ideals of
the matrix $\Gamma_{n}-$semiring $S_{n}$.
\label{Th:3.19}\end{Theorem}

\textit{Proof.} Let L be the left operator semiring of S. Then by
Theorem 5.2 of \cite{re:Biswas} there is an inclusion preserving
bijection between FI(L) and $FI(L_{n})$. In view of Theorem
\ref{Th:3.8}, there is an inclusion preserving bijection between
$FI(S_{n})$ and $FI([S_{n},\Gamma_{n}])$. Again by Proposition
\ref{Proposition:2.10}, $FI([S_{n},\Gamma_{n}])$ and $FI(L_{n})$ are
isomorphic. Further, there is an inclusion preserving bijection
between FI(S) and FI(L)\textit{(cf. Theorem \ref{Th:3.8})}.
Combining all these we obtain the theorem.

\textbf{Remark.} The above theorem can also be obtained directly via
the mapping $\mu \mapsto \mu_{n}$ where $\mu \in FI(S)$ and $\mu_{n}
\in FI(S_{n})$ and $\mu_{n}$ is defined by
\\$\mu_{n}([a_{ij}])=\min[\mu(a_{ij}):a_{ij} \in [a_{ij}]; 1\leq i,j
\leq n]$.

\begin{Remark} Throughout the paper we defined the concepts and
obtained results for left operator semiring. Similar concepts and
results can be obtained for right operator semiring R of S by using
analogy between L and R.\label{Re:3.12}\end{Remark}

\begin{Remark} From Theorem \ref{Th:3.8} and its dual (cf. Remark \ref{Re:3.12})
for right operator semiring R we easily obtain lattice isomorphism
of all fuzzy ideals of R and that of L.\end{Remark}
\textbf{Concluding Remark.} It is well known that operator semirings
of a
 $\Gamma-$semiring  are very effective in the study of
$\Gamma-$semirings. Lemma \ref{Lemma:3.11}, Lemma \ref{Lemma:3.12}, Theorems \ref{Th:3.15} and \ref{Th:3.19}
illustrate that operator semirings of a $\Gamma-$semiring can also
be made effective in the study of $\Gamma-$semiring in terms of
fuzzy subsets.


\begin{thebibliography}{99}
\bibitem{re:Biswas1} Dutta, T. K. and Biswas, B. K.: \textit{Fuzzy Prime Ideals Of A
Semiring}; Bull. Malaysian Math. Soc.(Second Series) \textbf{17}
(1994), 9-16.

\bibitem{re:Biswas} Dutta, T. K. and Biswas, B. K.: \textit{Structure of fuzzy ideals of
semirings}; Bull. Calcutta Math. Soc. \textbf{89(4)} (1997),
271-284.

\bibitem{re:Dutta} Dutta, T.K. and Sardar, S.K.: \textit{On the Operator Semirings of a
$\Gamma-$semiring}; Southeast Asian Bull. of Math.
\textbf{26}(2002), 203-213.

\bibitem{re:Chanda} Dutta, T.K. and Chanda, T.: \textit{Structures of Fuzzy Ideals of
$\Gamma-$Ring}; Bull. Malays. Math. Sci. Soc. (2) \textbf{28(1)}
(2005), 9-18.
\bibitem{re:Matrix} Dutta, T.K. and Sardar, S.K.: \textit{On Matrix $\Gamma-$semirings};
Far East J.Math. Sci.(FJMS) \textbf{7(1)} (2002), 17-31.

\bibitem{re:Goswami} Dutta, T.K.,  Sardar, S.K. and Goswami, S.: \textit{An introduction to fuzzy ideals of
$\Gamma-$semirings}; (To appear) Proceedings of \textit{National
Seminar on Algebra, Analysis and Discrete Mathematics.}, University of Kerala, India.

\bibitem{re:Goswami2} Dutta, T.K.,  Sardar, S.K. and Goswami, S.: \textit{Operations on fuzzy ideals of
$\Gamma-$semirings}; Communicated.

\bibitem{re:Golan} Golan, J.S.: \textit{Semirings and their applications}, Kluwer
Academic Publishers,1999.

\bibitem{re:Rao} Rao, M.M.K.: \textit{$\Gamma-$semiring-1}; Southeast Asian
Bull. of Math. \textbf{19} (1995), 49-54.

\bibitem{re:Sardar} Sardar, S.K.: \textit{On $\Gamma-$semifields};
Communicated.

\bibitem{re:Zadeh} Zadeh, L.A.:  \textit{Fuzzy sets}; Information and Control
\textbf{8}( 1965 ), 338-353.

\end{thebibliography}
\end{document}